\documentclass{amsart}
\usepackage{amssymb}
\newcommand{\fto}{{\mathsf {FtO}}}
\newcommand{\nor}[2]{\left|\!\left|#2\right|\!\right|_{#1}} 
\newcommand{\pnor}{\nor}

\newcommand{\Pg}{{\mathbb P}}

\newcommand{\forces}{\Vdash}
\newcommand{\ZFCa}{{\operatorname{\mathsf {ZFC}}}}
\newcommand{\CH}{\operatorname{\mathsf {CH}}}

\newcommand{\GCH}{\operatorname{\mathsf {GCH}}}
\newcommand{\NCF}{\operatorname{\mathsf {NCF}}}

\newcommand{\rest}{{\mathord{\restriction}}}

\newcommand{\Pp}{{\mathcal P}}
\newcommand{\HH}{{\mathbf H}}

\newcommand{\V}{{\mathbf V}}

\newcommand{\QED}{\hspace{0.1in} \square \vspace{0.1in}}

\newcommand{\lft}[2]{\mathopen\ifcase#1{}\oo\or
                        \big#2\or\Big#2\else\oo\fi} 
\newcommand{\rgt}[2]{\mathclose\ifcase#1{}\oo\or
                        \big#2\or\Big#2\else\oo\fi}

\newcommand{\C}{{\mathbf   {C}}}

\theoremstyle{plain}
\newtheorem{theorem}{Theorem}
\theoremstyle{plain}
\newtheorem{lemma}[theorem]{Lemma}

\newtheorem{definition}[theorem]{Definition}

\begin{document}
\title{There may be no Hausdorff ultrafilters}
\author{Tomek Bartoszynski}
\address{Department of Mathematics and Computer Science\\
Boise State University\\
Boise, Idaho 83725 U.S.A.}
\thanks{This material is based on work supported while the first author was
serving at NSF, and by 
NSF grant DMS 0200671 and by KBN grant 5 P03A 037 20.} 
\email{tomek@math.boisestate.edu, http://math.boisestate.edu/\char 126 tomek}
\author{Saharon Shelah}
\thanks{The second author was partially supported by Israel Science
   Foundation. Publication 826}
\address{Department of Mathematics\\
Hebrew University\\
Jerusalem, Israel}
\email{shelah@math.huji.ac.il, http://math.rutgers.edu/\char 126 shelah/}
\keywords{Hausdorff, ultrafilters, consistency}
\subjclass{03}
\begin{abstract}
An ultrafilter $U$ is Hausdorff if for any two functions $f,g \in
\omega^\omega$, $f(U)=g(U)$ iff $f \rest X=g\rest X$ for some $X \in U$. 
We will show that it is consistent that there are no Hausdorff ultrafilters.
\end{abstract}
\maketitle

\section{Introduction}
For $f \in \omega^\omega$ and an ultrafilter $U$ on $\omega$
define  
$f(U)=\{X \subseteq \omega: f^{-1}(X) \in U\}$. 
Let $\fto$ be the collection of all finite-to-one functions $f \in
\omega^\omega$. 
\begin{definition}
  Let $U$ be an ultrafilter on $\omega$. We say that
  \begin{enumerate}
  \item $U$ is Hausdorff if for any two functions $f, g \in
\omega^\omega$, if $f(U)=g(U)$ then $f \rest X=g\rest X$ for some $X \in U$.
\item $U$ is weakly Hausdorff if for any two functions $f, g \in
\fto$, if $f(U)=g(U)$ then $f \rest X=g\rest X$ for some $X \in U$.
  \end{enumerate}
\end{definition}

It is easy to see that
\begin{enumerate}
\item Ramsey ultrafilters are Hausdorff.
\item $q$-points  are weakly Hausdorff.
\item Weakly Hausdorff $p$-points are Hausdorff.
\end{enumerate}

It is worth mentioning that the following appears as an exercise in
\cite{je:settheory}. 
\begin{lemma}
If $f(U)=U$ then there exists $X \in U$ such that $f(n)=n $ for $n \in X$.
\end{lemma}
Therefore, if $U$ is not Hausdorff, then this is witnessed by two functions,
both 
not one-to-one 
mod $U$.

The notion of a Hausdorff ultrafilters was reintroduced and studied by Mauro
Di Nasso, Marco Forti  and others in a sequence of papers
(\cite{nassoforti},
\cite{nassoforti2}, \cite{nassobenci02} and \cite{nassoufilters}) in context
of topological extensions. 
They used the name  Hausdorff because Hausdorff 
ultrafilters are precisely those ultrafilters whose
ultrapowers equipped with the standard topology are Hausdorff topological
spaces. They 
asked whether the existence of a Hausdorff ultrafilter can be proved in
$\ZFCa$. We will show that, at least for ultrafilters on $\omega$, the answer
is negative.
However such ultrafilters (with various extra properties) may be constructed
under  from additional set theoretical assumptions (see \cite{nassoforti2}).

\section{Construction of the model}

In this section we will show how to build a model where there are no
Hausdorff ultrafilters  modulo the proofs of theorems \ref{one} and
\ref{two} below.  

\begin{definition}
  An ultrafilter $U$ is strongly non-Hausdorff if for every $f \in
  \fto$,
$f(U)$ is not weakly Hausdorff.
\end{definition}

\begin{theorem}\label{one}
  Assume $\CH$. There exists a strongly non-Hausdorff $p$-point.
\end{theorem}

\begin{definition} \cite{Bla86Nea}, \cite{Bla87Nea}, \cite{BlaShel89Nea}.
  Let ${\NCF}$ stand for the following statement:

for any ultrafilters $U,V $ on $\omega$ there exists $h \in \fto$
such that $h(U)=h(V)$.
\end{definition}

\begin{theorem}\label{two}
  There exists a proper forcing notion $\Pp$ such that
  \begin{enumerate}
  \item If $\V \models \GCH$ then $\V^{\Pp} \models
  2^{\aleph_0}=\aleph_2$,
\item If $U$  is a $p$-point  then 
$\V^{\Pp} \models U$ generates a $p$-point.
\item If $U$ is strongly non-Hausdorff filter then 
$\V^{\Pp} \models U$ generates a  strongly non-Hausdorff
  filter.
\item $\V^{\Pp} \models \NCF$.

  \end{enumerate}
\end{theorem}

\begin{theorem}
  Suppose that $\V \models \GCH$. Then in $\V^{\Pp}$ there are no
  weakly Hausdorff ultrafilters. In particular, there are no Hausdorff
  ultrafilters in this model.
\end{theorem}
\begin{proof}
  Let $U_0$ be a strongly non-Hausdorff $p$-point in $\V$ given by
  theorem \ref{one}. By theorem \ref{two}, $U_0$ generates  a strongly
  non-Hausdorff $p$-point in $\V^{\Pp}$, and  $\V^{\Pp}$ satisfies
  $\NCF$.
So suppose that $U$ is an ultrafilter in $\V^{\Pp}$. 
By $\NCF$ there exists $h \in \fto$ such that $h(U)=h(U_0)$. Since $U_0$ is
  strongly non-Hausdorff in $\V^{\Pp}$ it follows that
  $h(U_0)$ is not Hausdorff.
On the other hand if $U$ was Hausdorff then the following lemma would imply
  that $h(U)$ is Hausdorff as well, a contradiction.
\begin{lemma}\label{easy}
  If $U$ is Hausdorff then $h(U)$ is also Hausdorff.
\end{lemma}
\begin{proof}
  Let $f,g \in \fto$ be such that $f(h(U))=g(h(U))$. It follows that there is
  $X\in U$ such that $f \circ h\rest X=g\circ h\rest X$. Thus $f \rest h[X] =
  g \rest h[X]$ and $h[X] \in h(U)$.
\end{proof}
\end{proof}

\section{A strongly non-Hausdorff ultrafilter}


Let $I \subset \omega$ be a finite set and let $\Delta=\{(n,n) : n \in
\omega\}$.  Denote by $[I]^2=(I \times I) \setminus \Delta$. 
For a set $X \subseteq [I]^2$ define
$$\pnor{I}{X}=\min\left\{k: \exists \{A_i,B_i: i \leq k\} \ \forall i \leq k \
A_i \cap B_i =\emptyset \ \text{ and } X \subseteq \bigcup_{i\leq k}
A_i \times B_i\right\}.$$
We will drop the subscript $I$ if it is clear from the context what it
is.

\begin{lemma}\label{imp}
  \begin{enumerate}
  \item $\nor{I}{[I]^2} \longrightarrow \infty$ as $|I| \rightarrow \infty$.
  \item $\nor{I}{X \cup Y} \leq \nor{I}{X}+\nor{I}{Y}$,
  \item if $Z \subseteq I$ and $X \subseteq [I]^2$, $\nor{I}{X}>2$, then
  either $\nor{I}{[Z]^2\cap X} \geq \nor{I}{X}/2-1$
  or
$\nor{I}{[I \setminus Z]^2 \cap X} \geq
  \nor{I}{X}/2-1.$
  \end{enumerate}
\end{lemma}
\begin{proof}
  If (1) fails then there is $k \in \omega$ and sets $\{A^n_j, B^n_j: n,
  j\leq k\}$ such that $A^n_j\cap B^n_j=\emptyset$ for $j \leq k$ and
$[n]^2=\bigcup_{j\leq k} A_j^n \times B^n_j$. By compactness we get sets
$\{A^\omega_j, B^\omega_j: 
  j\leq k\}$ such that $A^\omega_j\cap B^\omega_j=\emptyset$ for $j \leq k$
and $[\omega]^2=\bigcup_{j\leq k} A_j^\omega \times B^\omega_j$, which is not
  possible. 

A more direct argument shows that
the following strategy is optimal for covering $[I]^2$, when
      $|I|$ is a power of two. Write 
      $I=I_0\cup I_1$ of equal size and use $I_0\times I_1$ and $I_1
      \times I_0$ to cover part of $I \times I$. For the rest, that is
      $(I_0 \times I_0) \cup (I_1 \times I_1)$
apply the same strategy by writing $I_0=I_{00} \cup I_{01}$ and
      $I_1=I_{10} \cup I_{11}$. The procedure terminates when squares
      have size $2 \times 2$, that is after $\log_2(|I|)-1$ steps. 
At that time we have used
$2+2\cdot 2 + 2 \cdot 4 + \cdots + 2 \times 2^{\log_2(|I|)-1}=2\cdot
      |I|-2$ rectangles. For $I$ of arbitrary size we get (by rounding
      down to the nearest power of two) that 
$\pnor{I}{[I]^2}\geq |I|-2$.

(2)  is obvious.

(3) Note that
\begin{multline*}
\pnor{I}{X} \leq \pnor{I}{\left([Z]^2  \cup [I \setminus Z]^2
\cup (Z \times (I
  \setminus Z))  \cup ((I \setminus Z) \times Z)\right) \cap X}   \leq \\
\pnor{I}{[Z]^2\cap X} + 
\pnor{I}{[I \setminus Z]^2\cap X }
 + 1 +1.
\end{multline*}
Thus 
$$\pnor{I}{[Z]^2\cap X} + 
\pnor{I}{[I \setminus Z]^2\cap X}
\geq \pnor{I}{X}-2.$$
\end{proof}

For $I \in [\omega]^{<\omega}$ let $\pi_1, \pi_2: [I]^2 \longrightarrow I$ be
projections onto first and second coordinate respectively.
\begin{lemma}\label{A}
  Suppose that $X \subseteq [I]^2$, and $\nor{I}{X} > 2$. Then $\pi_0(X)
  \cap \pi_1(X) \neq \emptyset$.
\end{lemma}
\begin{proof}
  Suppose that $\pi_0(X)=u$ and $\pi_1(X)=v$. If $u \cap v=\emptyset$ then
$X \subseteq (u\times v) \cup (v \times u)$. Thus $\nor{I}{X}\leq 2$.
\end{proof}

Next we define functions $f^0, g^0 \in \fto$ that will
witness that 
ultrafilter $U_0$ that we are about to  construct is not Hausdorff.

Let $\{I_k,J_k: k \in \omega\}$ be two sequences of disjoint consecutive
intervals such that for $k \in \omega$,
\begin{enumerate}
\item $\nor{I_k}{[I_k]^2} \geq 2^{2^k}$,
\item $|J_k|=|[I_k]^2|$.
\end{enumerate}
Bijection implicit in (2) allows us to define projections $\pi^k_0, \pi^k_1:
J_k 
\longrightarrow I_k$. Let 
 $f^0=\bigcup_k \pi^k_0$ and $g^0=\bigcup_k \pi^k_1$. Note that $f^0(x)\neq
 g^0(x)$ for any $x \in J_k=[I_k]^2$, $k \in \omega$.

As a warm-up let us use these definitions to show the following:
\begin{lemma}\label{h}
  Assume $\CH$. There exists a $p$-point that is not weakly Hausdorff.
\end{lemma}
\begin{proof}
We will need the following easy observation:

\begin{lemma}\label{B}
  If $f,g \in \fto$ and $U$ is an ultrafilter then the following
  conditions are equivalent:
  \begin{enumerate}
  \item $f(U) \neq g(U)$,
  \item $f[X] \cap g[X] =\emptyset$ for some $X \in U$.~$\QED$
  \end{enumerate}
\end{lemma}

  We will build an ultrafilter $V_0$ on the set $\bigcup_k [I_k]^2$ which we
  identified with $\omega$.
Let $\{Z_\alpha: \alpha < \omega_1\}$ be enumeration of
$[\omega]^\omega$.

We will build by induction a sequence
$\{X_\alpha: \alpha < \omega_1\}$ so that
\begin{enumerate}
\item $\forall \beta< \alpha\ X_\alpha \subseteq^\star X_\beta$,
\item $X_{\alpha+1} \cap Z_\alpha=\emptyset$ or $X_{\alpha+1} \subseteq
  Z_\alpha$ 
  for all $\alpha$.
\item for every $\alpha< \omega_1$, $f^0[X_\alpha] \cap g^0[X_\alpha] \neq
  \emptyset$. 
\item for every $\alpha< \omega_1$, $\limsup_k \nor{I_k}{X_\alpha \cap
    J_k}=\infty$. 
\end{enumerate}
Let $V_0=\{X: \exists \alpha \ X_\alpha \subseteq^\star X\}$. 
Note that
the conditions (1) and (2) guarantee that $V_0$ is a
$p$-point, and lemma \ref{B} and  (3) implies that
$f^0(V_0)=g^0(V_0)$. 
Finally, (4) is the requirement that (by lemma \ref{A}) implies (3).

{\sc Successor step}. Suppose that $X_\alpha$ is given. Find a strictly
increasing sequence $\{\ell_k: k \in \omega\}$ such that the set 
$A=\{k: \nor{I_k}{X_\alpha \cap J_k}=\ell_k\}$ is infinite.
Let
$A_0=\{k: \nor{I_k}{X_\alpha \cap Z_\alpha \cap J_k}\geq \ell_k/2-1\}$ and
$A_1=\{k: \nor{I_k}{(X_\alpha \setminus Z_\alpha) \cap J_k}\geq \ell_k/2-1\}$.
By lemma \ref{imp}(3), one of these sets, say $A_0$, is infinite.
Let $X_{\alpha+1}=\bigcup_{k\in A_0} X_\alpha\cap Z_\alpha\cap J_k$. The other
case is the same.

{\sc Limit step}. Given $\{X_\beta: \beta<\alpha<\omega_1\}$ let
$\{\beta_k: k \in \omega\}$ be an increasing sequence cofinal in $\alpha$. 
By finite modifications we can assume that $X_{\beta_{k+1}} \subseteq
X_{\beta_k}$ for all 
$k$.
Build by recursion a strictly increasing sequence $\{u_k: k \in \omega\}$ such
that 
$$\forall k\ \forall j\leq k \ \exists i\in [u_k, u_{k+1}) \
  \nor{I_i}{X_{\beta_j}  \cap
  J_i} \geq k,$$
and let $$X_\alpha= \bigcup_k \left(X_{\beta_k} \cap \bigcup_{i \in [u_k,
  u_{k+1})} 
  J_i\right).$$ 
It is clear that $X_\alpha$ satisfies (1) and (4).
\end{proof}

Observe that $\CH$ was only needed in the limit step.
If we do not require that that $U$ is a $p$-point then we have the following:
\begin{theorem}
  There exists an ultrafilter that is not weakly Hausdorff.
\end{theorem}
\begin{proof}
  As in lemma \ref{h},  we will build an ultrafilter on the set $\bigcup_k
  [I_k]^2$. 
Let 
$${\mathcal I}=\left\{X \subseteq \bigcup_k
  [I_k]^2: \limsup_k \nor{I_k}{X \cap
    J_k}<\infty\right\}.$$
Note that ${\mathcal I}$ is an ideal, and let $U$ be any ultrafilter
  orthogonal to ${\mathcal I}$. Functions $f^0, g^0$ witness that $U$ is not
  Hausdorff. 
\end{proof}

{\sc Proof of Theorem \ref{one}}. 

Now we are ready to construct  a $p$-point ultrafilter $U_0$ whose all finite-to-one
images are not weakly Hausdorff. 

Let $\{h_\alpha,  Z_\alpha: \alpha < \omega_1\}$ be enumeration of
$\fto$ and $[\omega]^\omega$ respectively.
We will build by induction sequences $\{ f^\alpha, g^\alpha:
\alpha<\omega_1\}$, 
$\{X_\alpha: \alpha < \omega_1\}$ so that
\begin{enumerate}
\item $\forall \beta< \alpha\ X_\alpha \subseteq^\star X_\beta$
\item $X_{\alpha+1} \cap Z_\alpha=\emptyset$ or $X_{\alpha+1} \subseteq
  Z_\alpha$ 
  for all $\alpha$.
\item for every $\alpha$, $f^{\alpha+1}, g^{\alpha+1}$ witness that
  $h_\alpha(U_0)$ is not Hausdorff, where $U_0=\{X \subseteq \omega:
  \exists \alpha \ X_\alpha 
\subseteq^\star X\}$. 
\item $\forall \beta \ \forall \alpha \geq \beta\ f^{\beta+1}\circ
  h_\beta[X_\alpha] \cap  g^{\beta+1}\circ
  h_\beta[X_\alpha] \neq \emptyset$.
\end{enumerate}
As before, (1) and (2) guarantee that $U_0$ is a
$p$-point, and
(3) implies that $U_0$ is strongly non-Hausdorff, and (4) is a specific form of
(3). 
Note that at the limit stages we only have to
preserve the induction hypothesis.
At the successor step we will first define an auxiliary function
$e_{\alpha+1}$, and put 
$f^{\alpha+1}=f^0\circ e_{\alpha+1}:h_\alpha[X_\alpha]\longrightarrow \omega$
and 
$g^{\alpha+1}=g^0\circ e_{\alpha+1}:h_\alpha[X_\alpha]\longrightarrow \omega$.
In other words, $f^{\alpha+1}, g^{\alpha+1}$ are copies of $f^0, g^0$
on the image  of $X_{\alpha}$ via $e_{\alpha+1}\circ h_\alpha$. 

Therefore we need to clarify condition (3) by imposing conditions on
$e_\alpha$ and specifying the induction hypothesis. 


\begin{definition}\label{wit}
Let us say that a finite set $Y \subseteq X_\alpha$ is a
$(n,\beta,\alpha)$-witness 
if   there exists $k \in \omega$ such that
$\nor{I_k}{e_{\beta+1}\circ h_\beta[Y] \cap J_k}\geq n$.
\end{definition}

To satisfy (3),  we demand that for $\beta< \alpha<\omega_1$,
\begin{enumerate}
\item[(5)] 
$\limsup_k \nor{I_k}{e_{\beta+1}\circ h_\beta[X_\alpha] \cap J_k}=\infty,$ or equivalently
\item[(6)] $\forall n \ \exists Y \in [X_\alpha]^{<\omega} \ Y \text{ is a $(n,\beta,\alpha)$-witness}$.
\end{enumerate}

{\sc  Limit step}.

Suppose that $\{X_\beta: \beta<\alpha\}$ are defined and $\alpha$ is a limit
ordinal.
Let
$\{\beta_k: k \in \omega\}$ be an increasing sequence cofinal in $\alpha$, and
let $\{\gamma_k: k \in \omega\}$ be an enumeration of $\alpha$ such that
$\gamma_j \leq \beta_k$ for $j \leq k$. Without loss of generality we can
assume that $X_{\beta_n} \subseteq X_{\beta_m}$ for $n \geq m$.

Build by recursion a strictly increasing sequence $\{u_k: k \in \omega\}$ such
that 
$$\forall k\ \forall l,j\leq k \ \exists i\in [u_k, u_{k+1}) \
  \nor{I_i}{e_{\gamma_l+1}\circ h_{\gamma_l}[X_{\beta_j}] \cap
  J_i} \geq k,$$
and let
 $$X_\alpha= \bigcup_k \left(X_{\beta_k} \cap \bigcup_{i \in [u_k, u_{k+1})}
  J_i\right).$$ 
It is clear that $X_\alpha$ satisfies (1) and (4).

{\sc Successor step}.

Suppose that $X_\alpha$ satisfying (4) is already defined and we want to define
$X_{\alpha+1}$ and $e_{\alpha+1}$ satisfying (2) and (5). 
Recall that by the induction hypothesis, for $\beta< \alpha$,
$$\forall n \ \exists Y \in [X_\alpha]^{<\omega} \ Y \text{ is a $(n,\beta,\alpha)$-witness}.$$

Let $\{\beta_k: k \in \omega\}$ be an enumeration of $\alpha$.
Find a sequence $\{E_k: k \in \omega\}$ of consecutive intervals
such that
\begin{enumerate}
\item $\forall k \ \forall j\leq k \ h_\alpha^{-1}(E_k)$ contains a $(k,\beta_j,\alpha)$-witness.
\item $\forall k \ E_k \cap h_\alpha[X_\alpha] \neq \emptyset$.

\end{enumerate}

Let
$e_{\alpha+1}(j)=k \iff j\in E_k$ for $j \in \omega$. Condition (2) implies
that $e_{\alpha+1} \circ h_\alpha[X_\alpha]=\omega$. In particular,
either 
$\limsup_k \nor{I_k}{e_{\alpha+1}\circ h_\alpha[Z_\alpha \cap X_\alpha] \cap
  J_k}=\infty,$ or $\limsup_k \nor{I_k}{e_{\alpha+1}\circ h_\alpha[X_\alpha
  \setminus 
Z_\alpha] \cap 
  J_k}=\infty $.
Let $X_{\alpha+1}$ be the appropriate set.
It remains to  check that
for $\beta \leq \alpha$,
$\limsup_k \nor{I_k}{e_{\beta+1}\circ h_\beta[X_{\alpha+1}] \cap J_k}=\infty.$
If $\beta=\alpha$ it follows immediately from the definition of
$X_{\alpha+1}$, so suppose that $\beta=\beta_j<\alpha$.

Let $m \in e_{\alpha+1}\circ h_\alpha[X_{\alpha+1}]\setminus j$ and note that
$(e_{\alpha+1}\circ h_\alpha)^{-1}(m)$ contains a $(m,\beta_j,\alpha)$-witness.
It follows that
$\limsup_k \nor{I_k}{e_{\beta+1}\circ h_\beta[X_{\alpha+1}] \cap J_k}=\infty,$
which finishes the construction.

\section{Forcing}

Since known models for $\NCF$ are obtained by  countable support iteration
we will look for a proper forcing notion $\Pg$ such that the iteration of $\Pg$
has the required properties. 

Suppose that $\Pg$ is a proper forcing notion. $\Pg$ preserves non-meager
  sets if for every countable elementary submodel $N \prec \HH(\chi)$
  containing $\Pg$, a condition $p \in \Pg\in N$ and a Cohen real $c$ over $N$
  there exists $q \geq p$ such that $q$ is $(N,\Pg)$ generic and
$q \forces_{\Pg} c \text{ is Cohen over } N[\dot{G}].$

By \cite{BJbook}, 6.3.16 this is equivalent to the property
$\sqsubseteq^{\textup{Cohen}}$ defined in \cite{BJbook}.

Let $\Pg$ be a proper forcing notion such that:
\begin{enumerate}
\item $\Pg$ preserves $p$-points,
\item $\Pg$ preserves non-meager sets, that is  it preserves $\sqsubseteq^{\textup{Cohen}}$.
\end{enumerate}

Let $\Pp=\Pg_{\omega_2}$ be the countable support iteration of $\Pg$ of length
$\omega_2$. We have the following:
\begin{enumerate}
\item $\Pp$ preserves $p$-points (see \cite{BlaShel87The} or
  \cite{BJbook} 6.2.6), 
\item $\Pp$ preserves non-meager sets (\cite{BJbook}, 6.3.20).
\end{enumerate}

Recall that if $\Pg$ is either Blass-Shelah forcing from \cite{BlaShel87The}
or Miller superperfect forcing, then $\Pp$ has the above properties (7.3.46
and 7.3.48 of \cite{BJbook}) and
$\V^{\Pp} \models \NCF$, \cite{BlaShel89Nea} or \cite{BlaShel87The}.
Therefore the following theorem concludes the proof of theorem \ref{two}.

\begin{theorem}
  Suppose that $\Pp$ is a proper forcing that preserves non-meager sets and
  $U_0 \in \V$ is a strongly non-Hausdorff ultrafilter defined earlier.
Then 
$$\V^{\Pp} \models U_0 \text{ generates a a strongly non-Hausdorff 
  filter} .$$
\end{theorem}
\begin{proof}
  Clearly, $U_0$ may not generate an ultrafilter in the extension, for example
  when $\Pp$ is Cohen forcing.

Let $\C$ be the Cohen forcing interpreted as adding a 
 function $c \in \fto$. Specifically, the conditions are finite sequences of
 consecutive intervals $\{I_k: k<n\}$ and $c(i)=k \iff i \in I_k$.

\begin{lemma}\label{last}
Let $c$ be  Cohen reals over $\V$.
Then for every $h \in \V \cap \fto$,
$h(U_0)$ is not Hausdorff as witnessed by $f^h=f^0 \circ c$ and $g^h=g^0
\circ c$. 
\end{lemma}
\begin{proof}
  This is quite easy. Given $s=\{L_k: k<n\} \in \C$, $X=X_\alpha \in U_0$ we
  extend $s$ by adding an interval $L_n$ so large
  that $L_n \supseteq (e_{\alpha+1})^{-1}(k)$ for some $k>n$. That means that 
$f^0 \circ c$ and
$g^h=g^0
\circ c$ agree with $f^0 \circ e_{\alpha+1}$ and $g^0 \circ e_{\alpha+1}$ on
  long enough segments to witness that
$f^0 \circ c[X] \cap g^0 \circ c[X]\neq \emptyset.$
\end{proof}

Let $\dot{h}$ be a $\Pp$-name for an element of $\fto$.
Let $N \prec \HH(\chi)$ be a countable submodel containing $U_0, \dot{h}, p, \Pp$
and let $c \in \V \cap \fto$ be a Cohen real over $N$.
Since $\Pp$ preserves non-meager sets there is $q \geq p$ which is $(N,\Pp)$
generic and 
$q \forces_{\Pp} c \text{ is Cohen over } N[\dot{G}].$ In particular, by lemma
\ref{last}, 
$$q \forces_{\Pp} \dot{h}(U_0) \text{ is not Hausdorff as witnessed by $f^0 \circ c$ and $g^0
\circ c$. }$$
By elementarity, it means that
$\V^{\Pp} \models h(U_0) \text{ is not Hausdorff.}$.
\end{proof}

{\bf Acknowledgements}:
We would like to thank Mauro Di Nasso and Marco Forti for introducing us to
the problem, and many helpful comments concerning the draft versions of this
paper.


\begin{thebibliography}{10}

\bibitem{BJbook}
Tomek Bartoszynski and Haim Judah.
\newblock {\em Set {T}heory: on the structure of the real line}.
\newblock A.K. Peters, 1995.

\bibitem{Bla86Nea}
Andreas Blass.
\newblock Near coherence of filters. {I}. {C}ofinal equivalence of models of
  arithmetic.
\newblock {\em Notre Dame Journal of Formal Logic}, 27(4):579--591, 1986.

\bibitem{Bla87Nea}
Andreas Blass.
\newblock Near coherence of filters. {II}. {A}pplications to operator ideals,
  the {S}tone-{C}ech remainder of a half-line, order ideals of sequences, and
  slenderness of groups.
\newblock {\em Transactions of the American Mathematical Society},
  300(2):557--581, 1987.

\bibitem{BlaShel87The}
Andreas Blass and Saharon Shelah.
\newblock There may be simple {$P\sb {\aleph\sb 1}$}- and {$P\sb {\aleph\sb
  2}$}-points and the {R}udin-{K}eisler ordering may be downward directed.
\newblock {\em Annals of Pure and Applied Logic}, 33(3):213--243, 1987.

\bibitem{BlaShel89Nea}
Andreas Blass and Saharon Shelah.
\newblock Near coherence of filters. {III}. {A} simplified consistency proof.
\newblock {\em Notre Dame Journal of Formal Logic}, 30(4):530--538, 1989.

\bibitem{nassoufilters}
Mauro Di~Nasso.
\newblock Ultrafilter semirings and nonstandard submodels of the {S}tone-{C}ech
  compactification of the natural numbers.
\newblock submitted for the volume "Logic and its application in algebra and
  geometry", Contemporary Mathematics AMS, edited by Yi Zhang.

\bibitem{nassoforti2}
Mauro Di~Nasso and Marco Forti.
\newblock Hausorff ultrafilters.
\newblock to appear.

\bibitem{nassoforti}
Mauro Di~Nasso and Marco Forti.
\newblock Topological and {N}onstandard extensions.
\newblock to appear.

\bibitem{je:settheory}
Thomas Jech.
\newblock {\em Set Theory}.
\newblock Academic Press, 1978.

\bibitem{nassobenci02}
Viero~Benci Mauro Di~Nasso and Marco~Forti.
\newblock Hausdorff nonstandard extensions.
\newblock {\em Bol. Soc. Parana. Mat.}, 20, 2002.

\end{thebibliography}
 
\end{document}